\documentclass[10pt]{article}
\usepackage{amsmath}
\usepackage{amsfonts}
\usepackage{amssymb}
\usepackage{comment}
\usepackage{graphicx,color}
\long\def\proof#1{\removelastskip\vskip\baselineskip\relax\noindent{\it
Proof\if!#1!\else\ \ignorespaces#1\fi.\ }\ignorespaces}

\newcommand{\C}{{\mathbb C}}

\newcommand{\z}{\zeta}

\newcommand{\al}{\alpha}

\newcommand{\squareforqed}{\hbox{\rlap{$\sqcap$}$\sqcup$}}
\newcommand{\qed}{\ifmmode\squareforqed\else{\unskip\nobreak\hfil
\penalty50\hskip1em\null\nobreak\hfil\squareforqed
\parfillskip=0pt\finalhyphendemerits=0\endgraf}\fi}

\newcommand{\litem}{\par\noindent\dimen0=\parindent%
    \advance\dimen0 by-4pt
               \hangindent=\dimen0\ltextindent}

\newcommand{\ltextindent}[1]{\hbox to \hangindent{#1\hss}\ignorespaces}
\newcommand{\ltextjndent}[1]{\hbox to \hangindent{#1\hss}\ignorespaces\kern-1ex}

\begin{document}
\pagestyle{plain}

\title{Ap\'ery Acceleration of Continued Fractions}
\author{Henri Cohen\\ henri.cohen2@free.fr}

\maketitle
\begin{abstract}
  We explain in detail how to accelerate continued fractions (for constants
  as well as for functions) using the method used by R.~Ap\'ery in his proof
  of the irrationality of $\z(3)$. We show in particular that this can be
  applied to a large number of continued fractions which can be found in
  the literature, thus providing a large number of new continued fractions.
  As examples, we give a new continued fraction for $\log(2)$ and for
  $\z(3)$, as well as a simple proof of one due to Ramanujan.
\end{abstract}

\smallskip

\section{Introduction and Notation}

After R.~Ap\'ery announced that he had proved the irrationality of $\z(3)$,
the first explicit proof was found by D.~Zagier using suitable telescoping
series, in what was to become the celebrated Wilf--Zeilberger method (see
\cite{PWZ}). It was only a year later that Ap\'ery explained his
method, see for instance \cite{Bat-Oli} and \cite{Ape}, which was quite
different. Although considerable effort
was spent in trying to generalize his method to prove the irrationality of
other constants, the only small successes were first, the discovery by
G.~Rhin and the author \cite{Coh-Rhi} of a rapidly convergent continued
fraction for $\z(4)$, unfortunately insufficient to prove its irrationality
(which is of course well-known), and second, the proof by A.~Jeannin
\cite{Jea} of the irrationality of the sum of the inverses of the Fibonacci
numbers.

\smallskip

Nonetheless, an examination of the numerous continued fractions which can be
found in the literature, both for constants, but also for functions, it appears
that a large number can be accelerated using Ap\'ery's technique, giving
usually completely new and interesting continued fractions, but without any
Diophantine consequences. A particularly simple example is as follows.
The continued fraction
$$\log(1+z)=z/(1+1^2z/(2-z+2^2z/(3-2z+3^2/(4-3z+4^2/(5-4z+\cdots)))))\;,$$
valid for instance for $|z|<1$, is trivially obtained from the Taylor expansion
$\log(1+z)=z-z^2/2+z^3/3-\cdots$ by Euler's very classical transformation of
a series into a continued fraction, which converges at exactly the same speed
(essentially in $O(z^n)$). Applying Ap\'ery's method, we obtain
\emph{automatically} the new (but well-known) continued fraction
$$\log(1+z)=2z/(1(z+2)-1^2z^2/(3(z+2)-2^2z^2/(5(z+2)-3^2z^2/(7(z+2)-\cdots))))\;,$$
which is just as simple but converges much faster, essentially in
$O((z^2/(1+\sqrt{z+1})^4)^n)$ (for instance, for $z=1/2$ this is approximately
$O(98^{-n})$ instead of $O(2^{-n})$).

The aim of this paper is to explain this method
in great detail, together with a large number of examples. 

\medskip

To abbreviate, we write simply ``CF'' instead of ``continued fraction''.
In addition, we use the following notation for a CF $S$, which may
differ from notation used in other papers in the literature:
$$S=a(0)+b(0)/(a(1)+b(1)/(a(2)+b(2)/(a(3)+\cdots)))\;,$$
and we denote as usual by $p(n)/q(n)$ the $n$th convergent, so that
$p(0)/q(0)=a(0)/1$, $p(1)/q(1)=(a(0)a(1)+b(0))/a(1)$, etc..., and
for $u=p$ or $u=q$ we have $u(n+1)=a(n+1)u(n)+b(n)u(n-1)$.

When $a(n)$ and $b(n)$ are polynomials $A(n)$ and $B(n)$ for $n\ge n_0$,
we will write the continued fraction as
$$S=((a(0),a(1),\dotsc,a(n_0-1),A(n)),(b(0),b(1),\dotsc,b(n_0-1),B(n)))\;.$$
For instance, the last continued fraction above will simply be written as
$$\log(1+z)=((0,(2n-1)(z+2)),(2z,-n^2z^2))\;.$$
We will consider exclusively CFs of that type.

\section{Bauer--Muir Acceleration}

Contrary to general CFs, when $a(n)$ and $b(n)$ are polynomials for $n$
large, it is quite simple to determine the speed of convergence, at least
heuristically, see \cite{Bel-Coh} and \cite{Coh}. When the convergence is
fast, there is no need to improve the speed of convergence. However, when
the convergence is slow, for instance when $S-p(n)/q(n)=O(n^{-\al})$, it
is useful if not necessary to do so. The classical idea, theorized in
particular by Bauer and Muir, is to use a so-called \emph{modification} of the
$n$th \emph{tail} $\rho(n)$ of the CF, defined by
$\rho(n)=b(n)/(a(n+1)+b(n+1)/(a(n+2)+\cdots))$, so that
\begin{align*}S&=a(0)+b(0)/(a(1)+b(1)/(a(2)+\cdots+b(n-1)/(a(n)+\rho(n))))\;,\text{\quad and}\\
  S&=\dfrac{p(n)+\rho(n)p(n-1)}{q(n)+\rho(n)q(n-1)}\;.\end{align*}
Thus if we can find a function $r(n)$ which approximates $\rho(n)$, the
quantity $(p(n)+r(n)p(n-1))/(q(n)+r(n)q(n-1))$ will be a better approximation
to $S$ than $p(n)/q(n)$.

This observation leads to the two fundamental formulas of Bauer--Muir
acceleration. First, $\rho(n)$ trivially satisfies the recursion
$\rho(n)=b(n)/(a(n+1)+\rho(n+1))$, in other words
$\delta(n):=\rho(n)(a(n+1)+\rho(n+1))-b(n)=0$.
Thus, we want a function $r(n)$ such that $$d(n):=r(n)(a(n+1)+r(n+1))-b(n)$$
is as small as possible.

Second, it is immediate to create a new CF whose convergents are $p'(n)/q'(n)$
with $u'(n)=u(n)+r(n)u(n-1)$ for $u=p$ and $u=q$: after a small computation,
we find that such a CF is given by $(a'(n),b'(n))$, with
\begin{align*}
  a'(0)&=a(0)+r(0),\quad b'(0)=-d(0),\quad a'(1)=a(1)+r(1)\;,\\
  a'(n)&=a(n)+r(n)-r(n-2)d(n-1)/d(n-2)\quad\text{for $n\ge2$, and}\\
  b'(n)&=b(n-1)d(n)/d(n-1)\quad\text{for $n\ge1$}\;.
\end{align*}

It is also very easy to give a recipe for the choice of $r(n)$, depending
on the different speeds of convergence of the CF. Let us give three examples.

\smallskip

{\bf Example 1.} We start from the standard series
$\pi/4=1-1/3+1/5-1/7+\cdots$, which is immediately transformed by Euler
into the CF $$\pi/4=1/(1+1^2/(2+3^2/(2+5^2/(2+7^2/(2+\cdots)))))\;,$$
which using our notation is written $\pi/4=((0,1,2),(1,(2n-1)^2))$.
This CF of course converges like the series, i.e., $\pi/4-p(n)/q(n)=O(1/n)$.
We now look for $r(n)$ such that $d(n)$ is small for $n$ large, and we
immediately find that there are two polynomial solutions $r(n)=2n-3$ (giving
$d(n)=-4$), and $r(n)=-(2n-1)$ (giving $d(n)=0$, which seems ideal!).
Remembering that $r(n)$ is supposed to approximate
$\rho(n)=(2n-1)^2/(2+(2n+1)^2/(2+\cdots))$, the only plausible solution
is $r(n)=2n-3$. We thus find the new CF
$\pi/4=((1,4,5,6),(-1,4,(2n-1)^2))$, which satisfies
$\pi/4-p(n)/q(n)=O(1/n^3)$, so we indeed have accelerated our CF.

It is immediate to check that we can continue in this way, and obtain
CFs with $a(n)=4k+2$, $b(n)=(2n-1)^2$ for $n$ large, converging in
$O(1/n^{2k+1})$. This will be essential in Ap\'ery's method.

\smallskip

{\bf Example 2.} Here, we start from the standard series
$\log(2)=-\log(1-1/2)=(1/2)+(1/2)^2/2+(1/2)^3/3+\cdots$, which is immediately
transformed by Euler into the CF $\log(2)=((0,3n-1),(1,-2n^2))$, in other words
$$\log(2)=1/(2-2.1^2/(5-2.2^2/(8-2.3^2/(11-\cdots))))\;.$$
This converges like the series: $\log(2)-p(n)/q(n)=O(1/(n2^n))$, so already
relatively fast. Choosing here $r(n)=1-n$ gives the new CF
$$\log(2)=((1/2,3n),(1/2,-2n^2))\;,$$ which satisfies
$\log(2)-p(n)/q(n)=O(1/n^32^n)$, so a small polynomial gain. This does not
seem very interesting since the initial CF already converges reasonably
fast, but we will see that combined with Ap\'ery's method leads to a
series which converges \emph{much} faster, roughly in $1/34^n$.

\smallskip

{\bf Example 3.} Set $S=\int_0^\infty e^{-t}/(1+t)\,dt$. We have the
well-known CF $S=1/(2-1^2/(4-2^2/(6-3^2/(8-\cdots))))$, in other words
$S=((2n),(1,-n^2))$, with speed of convergence $S-p(n)/q(n)=O(e^{-4\sqrt{n}})$.
Using $r(n)=-n$, we find the new CF $S=((2n+1),(-1,-n(n+1)))$, in other
words $S=1-1/(3-1.2/(5-2.3/(7-3.4/(9-4.5/(11-\cdots)))))$. This has
exactly the same speed of convergence, so no acceleration, but Bauer--Muir
has given us a \emph{different} and interesting CF for the same quantity.

\smallskip

{\bf Important Remark.} To keep things simple (and to be able to reasonably
iterate Bauer--Muir in Ap\'ery's method), we must restrict to polynomials
or rational functions with \emph{rational} coefficients. It is easy to check
that because of this restriction, we can apply Bauer--Muir when the
convergence is in $O(1/n^{\al})$ (since $\al$ will always be rational), or
in $O(1/(B^nn^{\al}))$ \emph{when $B$ is a rational number}, and some other
types such as Example 3 above. In particular, it is \emph{not} applicable if
the convergence is in $O(1/B^n)$ with $B$ irrational.

\section{Ap\'ery's method}

We have seen above that for certain CFs it is possible to iterate the
Bauer--Muir process and still keep reasonably simple formulas.
Ap\'ery's method consists simply in combining this with a \emph{diagonal}
process. More precisely, it proceeds as follows. Let $(a(n),b(n))$ be
a CF with partial quotients $(p(n),q(n))$. We define by abuse of notation
$a(n,0)=a(n)$, $b(n,0)=b(n)$, $p(n,0)=p(n)$, and $q(n,0)$, and assume that
at step $l$ we have a CF $(a(n,l),b(n,l))$ with partial quotients
$(p(n,l),q(n,l))$. We use Bauer--Muir with a sequence $r(n,l)$ chosen such
that $d(n,l)=r(n,l)(r(n+1,l)+a(n+1,l))-b(n,l)$ is ``small'', and we thus
obtain a new CF $(a'(n,l),b'(n,l))$ with partial quotients
$u'(n,l)=r(n,l)u(n,l)+u(n-1,l)$ for $u=p$ and $u=q$. We do \emph{not} set
$u(n,l+1)=u'(n,l)$, but we shift $n$ by $1$, so we define
$u(n,l+1)=u'(n+1,l)=r(n+1,l)u(n+1,l)+u(n,l)$ for $u=p$ and $u=q$, hence
$a(n,l+1)=a'(n+1,l)$ and $b(n,l+1)=b'(n+1,l)$. Note that this shift by 1
is crucial. By induction this defines the 2-dimensional arrays
$(a(n,l),b(n,l),r(n,l),d(n,l),p(n,l),q(n,l))$.

Setting as usual $u=p$ or $u=q$, we summarize all the definitions and
recursions as follows, where to simplify notation, it is convenient to define
$R(n,l)=a(n+1,l)+r(n+1,l)$:
\begin{align*}
  u(n+1,l)&=a(n+1,l)u(n,l)+b(n,l)u(n-1,l)\;,\\
  u(n,l+1)&=r(n+1,l)u(n+1,l)+u(n,l)\;,\\
  R(n,l)&=a(n+1,l)+r(n+1,l)\;,\quad d(n,l)=r(n,l)R(n,l)-b(n,l)\;,\\
  a(n,l+1)&=R(n,l)-r(n-1,l)d(n,l)/d(n-1,l)\;,\\
  b(n,l+1)&=b(n,l)d(n+1,l)/d(n,l)\;.\end{align*}

The game now consists in \emph{walking} in a regular or semi-regular manner
in the two-dimensional arrays $u(n,l)$. The first idea is a \emph{diagonal}
walk, i.e., find the CF corresponding to $p(n,n)/q(n,n)$ or more generally
$p(n+m)/q(n+m)$ for small fixed $m$. This is easily seen to involve annoying
denominators, which in general give complicated formulas, but in the case
where these denominators cancel, these CFs are ideal and the most interesting.
To avoid denominators, we simply use a \emph{staircase} walk, in other
words $p(n,n)/q(n,n)$, $p(n+1,n)/q(n+1,n)$, $p(n+1,n+1)/q(n+1,n+1)$, etc...
(or more generally with the first argument $n$ and $n+1$ replaced by $n+m$ or
$n+m+1$), thus giving CFs with $(a(2n),b(2n))$ and $(a(2n+1),b(2n+1))$
given by different formulas, what I call \emph{period 2} CFs. In addition,
by \emph{contraction} it is immediate to check from the staircase CF whether
the diagonal one has simple coefficients. I call these CFs the Ap\'ery
accelerates of the initial one, and are usually the ones which converge
the fastest. The precise formulas are as follows, where as usual $u$ denotes
$p$ or $q$:
\begin{align*}
  u(n+m+1,n+1)&=R(n+1+m,n)u(n+m,n+1)\\
  &\phantom{=}-d(n+m+1,n)u(n+m,n)\;,\quad\text{and}\\
  u(n+m,n+1)&=R(n+m,n)u(n+m,n)+b(n+m,n)u(n+m-1,n)\;.\end{align*}

Finally, we can also do a \emph{vertical} walk, i.e., find the CF corresponding
to $p(0,n)/q(0,n)$ or more generally $p(m,n)/q(m,n)$ for small fixed $m$.
These will not converge fast, but provide completely new CFs, and I call them
the Ap\'ery \emph{duals} of the initial CF, since their corresponding
2-dimensional arrays will simply be the transpose of the initial ones.
The precise formula is
$$u(m,n+1)=(R(m+1,n-1)+r(m+1,n))u(m,n)-d(m+1,n-1)u(m,n-1)\;.$$

Note that in some cases the arrays are only defined for $n\ge l$ or $2n\ge l$
for instance, so in that case the Ap\'ery dual is not defined, although
the corresponding CF may converge, but usually not to (a M\"obius transform
of) the initial limit $S$.

\smallskip

An important remark must be made at this point. Since the initial terms
of a CF can obey different formulas than the generic term, the formulas
that are given above for by Ap\'ery's method are only valid for the
generic terms, i.e., give $(a(n),b(n))$ only for $n$ ``large'' (usually
in fact $n\ge2$). By the elementary theory of CFs, this means that
the limit of the Ap\'ery accelerate or of the Ap\'ery dual is of the form
$(AS+B)/(CS+D)$, where $S$ is the initial value, and $A$, $B$, $C$, and $D$
are small integers which can easily be found. One could avoid this by
precisely keeping track of each Bauer--Muir step, but this would be extremely
cumbersome. We will see this at work in all the examples below.

\section{Basic Examples}

\subsection{$\log(2)$ and $\psi(z+1/2)-\psi(z)$}

We begin by the slowly convergent series $\log(2)=1-1/2+1/3-1/4+\cdots$,
which by Euler is transformed into the CF $\log(2)=((0,1),(1,n^2))$.
Using Bauer--Muir, we find the following 2-dimensional arrays:
$$(a(n,l),b(n,l),r(n,l),d(n,l))=(2l+1,n^2,n-l-1,-(l+1)^2)\;.$$
After correcting for the initial terms, we find the staircase Ap\'ery
accelerate $\log(2)=((n),((1,1),(n^2,(n+1)^2)))$, using an evident notation
for period 2 CFs; luckily, the contraction of this is very simple and leads
to the CF
$$\log(2)=((0,3(2n-1)),(2,-n^2))=2/(3-1^2/(9-2^2/(15-3^2/(21-4^2/(27-\cdots)))))\;,$$ which converges like
$(1+\sqrt{2})^{-4n}$, so quite fast. The vertical Ap\'ery dual is
$((0,1),(1,n^2))$, in other words the initial CF which is therefore self-dual.

But we can do more without any additional computation: in the above formulas,
$n$ is a formal variable, so if we simply change $n$ into $n+z$ for $z\in\C$,
the same recursions are valid, so we can replace $n$ by $n+z$ in the
arrays. Now
$$\sum_{n\ge1}\dfrac{(-1)^{n-1}}{n+z}=\dfrac{\psi(z/2+1)-\psi(z/2+1/2)}{2}\;,$$
where $\psi(z)$ is the logarithmic derivative of the gamma function,
so the Euler transform of this gives the CF
$(\psi(z/2+1)-\psi(z/2+1/2))/2=((0,z+1,1),(1,(n+z)^2))$,
and using the formulas for Ap\'ery we deduce the period 2 CF
$(\psi(z/2+1)-\psi(z/2+1/2))/2=((0,2z+1,n+z),((1,1),((n+z)^2,(n+1)^2)))$
and after trivial modifications
$$\psi(z+1/2)-\psi(z)=((0,4z-1,n+2z-1),((2,1),((n+2z-1)^2,(n+1)^2)))\;.$$
Contrary to the case $z=0$, the contraction of this CF is not simple.
However, in the special case $z=3/4$, $b(2n)=((2n/2)+1/2)^2$ and
$b(2n+1)=(((2n+1)/2)+1/2)^2$ are given by the same formula, so this
simplifies to $\psi(5/4)-\psi(3/4)=((0,2,n+1/2),(2,(n+1)^2/4))$,
and since $\psi(5/4)-\psi(3/4)=4-\pi$, after an evident simplification
we obtain the following CF for $\pi$:
$$\pi=((4,4,2n+1),(-4,(n+1)^2))=4-4/(4+2^2/(5+3^2/(7+4^2/(9+5^2/(11+\cdots)))))\;,$$
which converges in $O((1+\sqrt{2})^{-2n})$.

In addition, we find that for $z\ne0$ the initial CF is not self-dual,
so after a small computation we find that the dual is still another CF:
$$\psi(z+1/2)-\psi(z)=((0,4z-1),(2,n^2))\;,$$
which converges in $O((-1)^n/n^{4z-1})$, to be compared with our initial
CF which converged in $O((-1)^n/n)$.

\smallskip

In this basic example, we see how using Ap\'ery, a single series can lead to
many interesting CFs.

\subsection{$\z(2)=\pi^2/6$ and $\psi'(z)$}\label{sec:zeta2}

Since this is similar, we will be brief. Starting from
$\z(2)=\pi^2/6=1+1/2^2+1/3^2+\cdots$, by Euler we have the trivial CF
$\z(2)=((0,2n^2-2n+1),(1,-n^4))$. Using Bauer--Muir, we obtain the arrays
\begin{align*}(a(n,l),b(n,l))&=(2n^2-2n+1+l^2+l,-n^4)\;,\quad\text{and}\\
  (r(n,l),d(n,l))&=(-n^2+(l+1)n-(l+1)^2/2,-(l+1)^4/4)\end{align*}
The Ap\'ery accelerate can be nicely contracted, and leads to Ap\'ery's
famous CF $$\z(2)=((0,11n^2-11n+3),(5,n^4))\;,$$ which converges in
$O((-1)^n((1+\sqrt{5})/2)^{-10n})$, and which can be used to prove
the (known) irrationality of $\z(2)$. The Ap\'ery dual is the CF
$((0,2n-1),(2,n^4))$, which corresponds to the alternating series
$\z(2)=2(1-1/2^2+1/3^2-\cdots)$.

Using the same recursions, one can also Ap\'ery accelerate the more general
series $\sum_{n\ge1}1/(n+z)^2=\psi'(z+1)$, and obtain in this way interesting
CFs for $\psi'(z)$.

We have seen that the alternating sum $1-1/2^2+1/3^2-\cdots$ is the dual
of the sum with positive signs, so its Ap\'ery accelerate is the same.
The corresponding function $\sum_{n\ge1}(-1)^{n-1}/(n+z)^2$ gives interesting
CFs for $\psi'(z/2)/2-\psi'(z)$.

\subsection{$\z(3)$ and $\psi''(z)$}

Starting from $\z(3)=1+1/2^3+1/3^3+\cdots$, by Euler we obtain the CF
$\z(3)=((0,(2n-1)(n^2-n+1)),(1,-n^6))$. Using Bauer--Muir, we find the
arrays
\begin{align*}(a(n,l),b(n,l))&=(2n^3-3n^2+(4l^2+4l+3)n-(2l^2+2l+1),-n^6)\;,\quad\text{and}\\
  (r(n,l),d(n,l))&=(-n^3+2(l+1)n^2-2(l+1)^2n+(l+1)^3,(l+1)^6)\;.\end{align*}
The Ap\'ery accelerate can be nicely contracted, and leads to Ap\'ery's
famous CF $$\z(3)=((0,(2n-1)(17n^2-17n+5)),(6,-n^6))\;,$$
which converges in $O((1+\sqrt{2})^{-8n})$ and which he used
together with Diophantine results to prove the irrationality of $\z(3)$.
The initial CF is self-dual.

The corresponding function $\sum_{n\ge1}1/(n+z)^3$ gives interesting
CFs for $\psi''(z)$. In particular, its Ap\'ery dual gives the CF
$$\z(3,z+1)=-\psi''(z+1)/2=((0,n^3+(n-1)^3+2z(z+1)(2n-1)),(1,-n^6))\;,$$
with speed of convergence $O(1/n^{4z+2})$, a formula due to Ramanujan.
This is possibly one of the simplest proofs of this formula.

\smallskip

\subsection{An Example where Ap\'ery can be Iterated}

\smallskip

We have mentioned that Bauer--Muir cannot reasonably be applied iteratively
when the convergence is in $O(1/B^n)$ with $B$ irrational. In particular,
in the three examples above the Ap\'ery accelerates have $B=(1+\sqrt{2})^4$,
$B=((1+\sqrt{5})/2)^{10}$, and $B=(1+\sqrt{2})^8$, so there is no hope
of applying Bauer--Muir and a fortiori Ap\'ery once again. However, there
are cases where it can be done, as in the following example.

We start with the following CF: $2^{1/3}=((1/2,7n-5),(1,-4n(3n-2)))$, which
converges in $1/((4/3)^nn^{5/3})$, so already exponentially fast with
$B=4/3$. Using the Ap\'ery formulas we obtain the period 2 CF
$((4n,4n+2),(-3n^2+2n,-3n^2-8n-5))$, where we do not give the initial terms
which have to be determined at the end (the complete CF is in fact
$(((1/2,2),(4n,4n+2)),((1/2,-5),(-3n^2+2n,-3n^2-8n-5)))$), but since we are
not finished it is a waste of time to do it here. This converges in
$1/(-3)^n$, which is already faster. To be able to continue, we need a period
1 CF, so we contract the above CF, and we are in luck, the simplifications is
still very simple: $((5(2n-1)),(-(9n^2-4)))$, which of course converges in
$O(1/(-3)^{2n})=O(1/9^n)$ (the complete CF is in fact $((1/2,3,5(2n-1)),(2,-(9n^2-4)))$).

Since $9$ is rational, we can try to use Ap\'ery once again, and indeed
it works, and another miracle happens, the period 2 CF does not need to
be contracted, it is already in fact a period 1 CF
$((9(2n-1)),(-(9n^2-16)))$, which converges in $O((-(1+\sqrt{2})^4)^{-n})$,
and since $B=-(1+\sqrt{2})^4$ is irrational, we cannot use Ap\'ery anymore.

To finish and obtain the complete CF, we need to find the initial terms of the
CF. To do this, we
proceed as follows. We start from the (almost certainly wrong) period 1
CF obtained above, and compute its limit $L$ \emph{numerically}. Since the
convergence is very fast, this is easy, and we find $L=-7.2728\cdots$.
We know that $2^{1/3}=(AL+B)/(CL+D)$ for (hopefully small) integers $A$, $B$,
$C$, $D$: to find them, we use a linear dependence algorithm such as LLL
on the numbers $(1,L,2^{1/3},2^{1/3}L)$, and we find immediately that
$2^{1/3}=(-L-13)/(2L+10)$, at least numerically. Note that, even though this
is a numerical hence non-rigorous check, it would be easy but tedious to
\emph{prove} it. Now the CF begins as $-9+16/(9+7/27+\cdots)$, so we
replace $L$ by $-9+16/(9+x)$, where $x$ represents the tail of the CF
starting at $7/(27+\cdots)$, and we find $2^{1/3}=(x+13)/(2x+10)=1/2+4/(x+5)$.
We thus create the CF $((1/2,5,9(2n-1)),(4,-(9n^2-16)))$, and we check that
its limit is indeed $2^{1/3}$, so
$$2^{1/3}=((1/2,5,9(2n-1)),(4,-(9n^2-16)))\;.$$

\section{Generalizations of Ap\'ery's Method}

\subsection{Case when $r(n,l)$ is a rational function in $l$}

Even though we look by indeterminate coefficients or otherwise for $r(n,l)$ as
a polynomial in $n$, it may happen that as a function of $l$ it is not a
polynomial (as in all the above examples), but a rational function of $l$.
Consider for instance the following (admittedly already rather complicated)
example. Set
$$G_3=L(\chi_{-3},2)=\sum_{n\ge1}\dfrac{(-3/n)}{n^2}=1-\dfrac{1}{2^2}+\dfrac{1}{4^2}-\dfrac{1}{5^2}+\cdots$$
It is possible to prove that we have the following CF:
$$G_3=((3/4,(2n-1)(9n^2-9n+22)),(2/3,-9n^4(9n^2-1)))\;,$$
which converges in $O(1/4^n)$. Using Bauer--Muir, we find the arrays
\begin{align*}(a(n,l),b(n,l))&=((2n-1)(9n^2-9n+18(l+1)^2+4),-9n^4(9n^2-1))\;,\\
  r(n,l)&=-9n^3+(18l+27)n^2-(18l^2+54l+40)n\\
  &\phantom{=}+(2/9)(3l+4)^2(3l+5)^2/(2l+3)\;,\text{\quad and}\\
  d(n,l)&=(4/81)(3l+4)^4(3l+5)^4/(2l+3)^2\;.\end{align*}
We see that $r(n,l)$ (hence also $d(n,l)$) is a rational function in $l$.
Note that the Ap\'ery accelerate and Ap\'ery dual, while completely
explicit, are rather complicated, but the main purpose of this example was to
show that $r(n,l)$ can be a rational function in $l$.

\subsection{Need for a simplifier}

Consider the following CF for $\pi$, which converges as $O(1/(4^nn^{3/2}))$:
$$\pi=((3,24,20n^2+4n+1),(3,-8n(2n+1)^3))\;.$$
Applying Bauer--Muir, we find that the new $a(n)$ has a denominator $2n-1$,
which indicates that if we continue in order to find our 2-dimensional arrays,
the degree of the denominator will increase, making the formulas unwieldy.
However, recall that if $(a(n),b(n))$ is a given CF, for any nonzero function
$t(n)$ with $t(0)=1$, the CF $(t(n)a(n),t(n)t(n+1)b(n))$ will have the same
convergents $p(n)/q(n)$ (more precisely, both $p(n)$ and $q(n)$ will be
multiplied by the same quantity $tf(n):=t!(n)=t(1)t(2)\cdots t(n)$). In our
particular case, choosing $t(n)=(2n-1)/(2n+1)$ not only gets rid of the
denominator, does not introduce new denominators, and also keeps constant the
degrees of $a(n)$ and $b(n)$. Thus, we generalize the Ap\'ery recursions
given above as follows:

\begin{align*}
  u(n+1,l)&=a(n+1,l)u(n,l)+b(n,l)u(n-1,l)\;,\\
  u(n,l+1)&=tf(n,l)(r(n+1,l)u(n+1,l)+u(n,l))\;,\\
  R(n,l)&=a(n+1,l)+r(n+1,l)\;,\quad d(n,l)=r(n,l)R(n,l)-b(n,l)\;,\\
  tf(n,l)&=t(n,l)tf(n-1,l)\;,\\
  a(n,l+1)&=t(n,l)(R(n,l)-r(n-1,l)d(n,l)/d(n-1,l))\;,\\
  b(n,l+1)&=t(n,l)t(n+1,l)b(n,l)d(n+1,l)/d(n,l)\;.\end{align*}

The formulas for the staircase walks are modified as follows:

\begin{align*}
  u(n+m+1,n+1)&=t(n+1+m,n)R(n+1+m,n)u(n+m,n+1)\\
  &\phantom{=}-tf(n+m+1,n)d(n+m+1,n)u(n+m,n)\;,\quad\text{and}\\
  u(n+m,n+1)&=tf(n+m,n)R(n+m,n)u(n+m,n)\\
  &\phantom{=}+tf(n+m,n)b(n+m,n)u(n+m-1,n)\;,\end{align*}
and the formula for the vertical walks is:
\begin{align*}u(m,n+1)&=(t(m+1,n-1)R(m+1,n-1)+r(m+1,n))u(m,n)\\
  &\phantom{=}-t(m+1,n-1)d(m+1,n-1)u(m,n-1)\;.\end{align*}

Coming back to our example of a CF for $\pi$, we find the arrays:
\begin{align*}
(a(n,l),b(n,l))&=(20n^2+(8l+4)n+1,-4(2n-l)(2n+1-l)(2n+1)^2)\;,\\
r(n,l)&=-4n^2+(8l+4)n-(12l^2+20l+9)\;,\\
d(n,l)&=-144(2n-l+1)(l+1)^3\;,\text{\quad and}\\
(t(n,l),tf(n,l))&=((2n-l-1)/(2n-l+1),-(l-1)/(2n-l+1))\;.\\
\end{align*}

Because of the denominator $2n-l+1$, the vertical walk gives a CF which
is unrelated to the initial one (and in fact converges to a limit related
to Catalan's constant). The staircase walks are OK since they us $l=n+m$,
which is smaller than $2n$ for $n$ large, and give semi-complicated period 2
CFs for $\pi$ converging like $(-1)^n/((1+\sqrt{5})/2)^{5n}$.

\subsection{Need for a multiplier: the Case of $\z(4)=\sum_{n\ge1}1/n^4$}

The next variant that we will study is perhaps the most interesting
(of course in addition to the basic one which is already extremely useful),
since it provides for instance a nice CF for $\z(4)$ (found in 1980 by
G.~Rhin and the author using exactly this method), and a new CF for $\z(3)$.

We can try to proceed as for $\z(2)$ and $\z(3)$. Applying Euler to the
series, we obtain the trivial CF $\z(4)=((0,2n^4-4n^3+6n^2-4n+1),(1,-n^8))$.
Applying Bauer--Muir creates a denominator $28(n-1)^2-9$ in $a(n)$ and
and $28n^2-9$ in $b(n)$. We could try to use a simplifier, but the numerator
of $a(n)$ being a degree 6 irreducible polynomial, this is hopeless.

Our salvation comes from an a priori bad idea of using a ``simplifier''
(here, the opposite, a ``complexifier'') on the initial CF, by multiplying
the CF by $2n-1$, in other words, using as initial series the more complicated
$$\z(4)=((0,(2n-1)(2n^4-4n^3+6n^2-4n+1)),(1,-(2n-1)(2n+1)n^8))\;.$$
The ``miracle'' is that when we now apply Bauer--Muir, the denominator
is simply $n-1$, and one easily checks that the simplifier $t(n)=(2n-2)/(2n)$
works. And it is immediate to check that this continues, so that we can
apply Ap\'ery with simplifier to this more complicated CF. We find:

\begin{align*}
a(n,l)&=(2n-1)\left(2n^4-4n^3+6n^2-4n+1+\dfrac{l(l+1)}{2}(17n^2-17n+5)\right)\;,\\
b(n,l)&=-n^6(4n^2-l^2)(4n^2-(l+1)^2)/4\;,\\
r(n,l)&=(2n+l)\left(-n^4+\dfrac{7}{2}(l+1)n^3-6(l+1)^2n^2+6(l+1)^3n-3(l+1)^4\right)\;,\\
d(n,l)&=-9(l+1)^8(4n^2-l^2)\;,\\
t(n,l)&=(2n-l-2)/(2n-l)\;,\\
tf(n,l)&=-l/(2n-l)\;.
\end{align*}

Once again since $l<2n$ the vertical walks are unrelated to $\z(4)$. On
the other hand, the staircase walks are well-defined, and for $m=0$ its
contraction is reasonably simple and gives the CF
$$\z(4)=\pi^4/90=((0,3(2n-1)(3n^2-3n+1)(15n^2-15n+4)),(13,3n^8(9n^2-1)))\;,$$
which converges like $O((2+\sqrt{3})^{-6n})$, not sufficiently fast to have
Diophantine applications, and as already mentioned, which
was discovered by G.~Rhin and the author in 1980, see \cite{Coh-Rhi}.

As usual, we can use the same recursions to find CFs for
$\psi'''(z)=6\sum_{n\ge0}1/(n+z)^4$ by replacing $z$ by $n+z-1$.
This gives rather complicated formulas. The only slightly interesting one
is the following, whose limit I do not know, but which gives the
\emph{asymptotic expansion} of $\psi'''(z)$ as $z\to+\infty$:
\begin{align*}
  \psi'''(z)&\sim((0,(2n-1)(n^4-2n^3-2(z^2-z-1)n^2+(2z^2-2z-1)n\\
  &\phantom{=}-(z-1)z(z^2-z+1))),(4z-2,-n^8(n^2-(2z-1)^2)))\;.\end{align*}
More precisely, changing $z$ into $1/z$ and simplifying by using $z^4$ as
simplifier, we find a CF such that $b(n)$ is divisible by $z^6$ for $n\ge1$,
so can be expanded in a power series, and we find
$$\psi'''(z)=\dfrac{2}{z^3}+\dfrac{3}{z^4}+\dfrac{2}{z^5}-\dfrac{1}{z^7}+\dfrac{4/3}{z^9}-\dfrac{3}{z^{11}}+\dfrac{10}{z^{13}}-\dfrac{691/15}{z^{15}}+\cdots\;,$$
which is indeed the asymptotic expansion of $\psi'''(z)$ for large $z$
(note that the coefficients of $1/z^{2k}$ vanish for $k\ne2$ as well as the
telltale presence of $691$, both sure signs of Bernoulli numbers).

\subsection{Need for a multiplier: a New CF for $\z(3)$}

We now give another example which leads to a new, but not especially
interesting, CF for $\z(3)$. Here we start with the alternating series
$\z(3)=(4/3)(1-1/2^3+1/3^3-1/4^3+\cdots)$, which by Euler is equivalent to the
CF $\z(3)=((0,3n^2-3n+1),(4/3,n^6))$. Once again, applying Bauer--Muir
directly leads to complicated denominators. Exactly as for the case of $\z(4)$,
we use a multiplier before applying Bauer--Muir.

Before giving the answer, I will explain how to find such factors without
guesswork. Let $n+u$ be the unknown linear factor. Using it as a multiplier,
we obtain the CF $((0,(n+u)(3n^2-3n+1)),(u+1,(n+u)(n+u+1)n^6))$. We now
apply the Bauer--Muir formulas, keeping $u$ as a formal variable. We
find for instance that $a(n)=N(n,u)/D(n,u)$ for complicated polynomials
$N$ and $D$ in the variables $n$ and $u$. Since we want this to simplify,
we compute the \emph{resultant} of $N(n,u)$ and $D(n,u)$ with respect to
the variable $n$. This gives a polynomial in $u$ alone, and since we want $u$
to be a rational number, we easily find that the rational roots of this
polynomial are $u=-1/2$ with multiplicity $4$, and $u=-39/14$ with multiplicity
$6$. Trying both, we find that $u=-39/14$ gives a denominator of degree 2,
probably difficult to simplify, while $u=-1/2$ gives a denominator $2n-3$,
which is promising.

Indeed, using as multiplier $2n-1=2(n-1/2)$, i.e., using the CF
$3\z(3)/4=((0,(2n-1)(3n^2-3n+1)),(1,(2n-1)(2n+1)n^6))$, we find the following
arrays:

\begin{align*}
(a(n,l),b(n,l))&=((2l+1)(2n-1)(3n^2-3n+1),n^6(4n^2-(2l+1)^2))\;,\\
r(n,l)&=(n-2(l+1))(2n+2l+1)(n^2-2(l+1)n+4(l+1)^2)\;,\\
d(n,l)&=-64(l+1)^6(4n^2-(2l+1)^2)\;,\text{\quad and}\\
(t(n,l),tf(n,l))&=((2n-(2l+3))/(2n-(2l+1)),-(2l+1)/(2n-(2l+1)))\;.\\
\end{align*}

Since $2n-2l+1$ never vanishes, the Ap\'ery dual is related to $\z(3)$,
and indeed we find that after simplification it gives the trivial CF
coming by Euler from the series with positive terms giving $\z(3)$, so
totally uninteresting. On the contrary, the Ap\'ery accelerate is a period
2 CF which can be reasonably simplified and gives the new CF:
$$\z(3)=((0,65n^4-130n^3+105n^2-40n+6),(7,-4(16n^2-1)n^6))\;,$$
which converges like $O(1/64^n)$, so quite fast, but not sufficiently fast
to have any Diophantine application.

Note that even though $64$ is a rational number, applying Bauer--Muir to this
CF leads to complicated denominators, so it is hopeless to try and use
Ap\'ery once again.

\subsection{Slowing Down Ap\'ery's Method}

The main goal of Ap\'ery's method is to accelerate as much as possible
a continued fraction. However, we have seen above that it is also very
useful for creating new CFs. Thus, we can try to use \emph{sub-optimal}
Bauer--Muir iterations, together with Ap\'ery's diagonal or staircase walks,
in order to find new CFs. We give a detailed example of how to proceed.

Consider the CF $\pi^2=((0,2n^2-2n+1),(6,-n^4))$, which is simply Euler's
transformation of the series $\pi^2=6\sum_{m\ge1}1/m^2$. Using Ap\'ery's
method using optimal Bauer--Muir leads immediately to his famous CF
for $\pi^2$, as we have seen in Section \ref{sec:zeta2}. The first
Bauer--Muir iteration is done with $r(n)=-n^2+n-1/2$, which is optimal
since it leads to a constant $d(n)$. Let us use instead $r(n)=-n^2+n+r_0$,
for some unknown $r_0$. To find a suitable $r_0$, we would like the
next $a(n)$ and $b(n)$ to remain polynomials. Keeping $r_0$ as a formal
variable, we thus compute $d(n)=r(n)(a(n+1)+r(n+1))-b(n)$, then
$a'(n)=a(n+1)+r(n+1)-r(n-1)d(n)/d(n-1)$ and $b'(n)=b(n)d(n+1)/d(n)$.

Both are rational functions in $n$ and $r_0$, and we want them to simplify,
so we again compute the \emph{resultant} with respect to $n$ of the numerator
and denominator of both, thus obtaining two polynomials in $r_0$, we compute
their GCD, and finally the roots, and we find that there are two possibilities
$r_0=-1/2$, which is the optimal one, and also $r_0=0$. We thus try
$r(n)=-n^2+n$, and we find $d(n)=n$ (which is of course not constant),
and $a'(n)=2n^2-n+1$, $b'(n)=-n^3(n+1)$. We check that the corresponding CF
has indeed been accelerated: it converges in $O(1/n^2)$ instead of $O(1/n)$
(the optimal one with $r_0=-1/2$ converges in $O(1/n^3)$).

Of course we now want to continue. Using the same method on this new CF,
we now find that we have three possibilities for $r_0$: $r_0=-2/3$, $0$,
and $2$. The value $r_0=-2/3$ corresponds to the optimal Bauer--Muir (i.e.,
with constant $d(n)$), but the other two values are plausible, and there is
no reason to choose one over the other. We decide (arbitrarily) to choose the
largest $r_0$, which will probably give the slowest acceleration, and in this
way we obtain very nice $2$-dimensional arrays:

\begin{align*}
(a(n,l),b(n,l))&=(2n^2+(l-2)n+1,-n^3(n+l))\;,\\
(r(n,l),d(n,l))&=(-n^2+n+l(l+1),(l+1)^3(n+l))\;.
\end{align*}

The vertical walk gives the initial CF, but the staircase walk can be nicely
contracted into a diagonal walk, and after correcting for the initial terms,
we find a new (but known) CF for $\pi^2$:
$$\pi^2=((0,5n^2-4n+1),(18,-2n^3(2n-1)))\;,$$
which converges in $O(1/(4^nn^{3/2}))$.

Since $B=4$ is rational, we can apply Ap\'ery's method once again
(slow or not), but we do not find any new CF but ``only'' Ap\'ery's CF.

\smallskip

We finish by giving another example of slow Ap\'ery, which leads to an
apparently completely new CF for $\log(2)$. Applying Euler's transformation
to the series $\log(2)=-\log(1-1/2)=(1/2)+(1/2)^2/2+(1/2)^3/3+\cdots$
gives the CF $\log(2)=((0,3n-1),(1,-2n^2))$, which of course converges
at exactly the same speed as the series, in $O(1/(n2^n))$. We proceed
exactly as before, in fact again with $r(n)=-n^2+n+r_0$ for an unknown $r_0$,
and we immediately find again very nice 2-dimensional arrays:
\begin{align*}
  (a(n,l),b(n,l))&=(3n^2-n+2l(2n-1),-2n^3(n+2l+1))\;,\\
  (r(n,l),d(n,l))&=(-n^2+n+2(l+1)(2l+1),4(n+2l+1)(n+2l+2)(l+1)^2)\;.
\end{align*}
The vertical walk gives the interesting but known CF
$$\log(2)=((0,4n^2-3n+1),(1/2,-2n^3(2n+1)))$$ which converges very slowly in
$O(1/n^{1/2})$. The ``miracle'' (which, as already mentioned, does not
occur very often) is that the staircase walk can again be nicely contracted
into a diagonal walk, and after correcting for the initial terms, we find
a new (and this time I believe completely new) CF for $\log(2)$:
$$\log(2)=((0,29n^2-29n+8),(5,-6n^2(9n^2-1)))\;,$$
which converges in $O(1/(27/2)^n)$, which is rather fast.

Since $27/2$ is rational, we can try to apply Ap\'ery's method once again
(slow or not), but here Bauer--Muir introduces denominators which cannot
be removed using simplifiers.

A similar procedure also gives the following new CF for $\log(2)$:
$$\log(2)=((0,59n^2-59n+20),(5,-24n^2(36n^2-1)))\;,$$
which converges in $O(1/(32/27)^n)$.


\bigskip

\begin{thebibliography}{14}

\bibitem{Ape}[Ape] R.~Ap\'ery, {\it Irrationalit\'e de $\z(2)$ et $\z(3)$\/},
  Ast\'erisque {\bf 61} (1979), 11--13.
\bibitem{Bat-Oli}[Bat-Oli] C.~Batut and M.~Olivier, {\it Sur
  l'acc\'el\'eration de la convergence de certaines fractions continues\/},
  S\'eminaire Th.~des Nombres Bordeaux (1979--1980), {\bf 9}, 1--26.
\bibitem{Bel-Coh}[Bel-Coh] K.~Belabas and H.~Cohen, {\it Numerical Algorithms
  for Number Theory using Pari/GP\/}, AMS Math.~Surveys and Monographs
  {\bf 254} (2021).
\bibitem{Coh}[Coh] H.~Cohen, {\it Continued fractions of polynomial type:
  theory and encyclopedic dictionary\/}, 420p., in preparation.
\bibitem{Coh-Rhi}[Coh-Rhi] H.~Cohen and G.~Rhin, {\it Acc\'el\'eration de la
  convergence de certaines r\'ecurrences lin\'eaires\/}, S\'eminaire de
  Th.~Nombres Bordeaux (1980/81), Expos\'e {\bf 16}.
\bibitem{Cuyt}[Cuyt] A.~Cuyt, V.~Petersen, B.~Verdonk, H.~Waadeland, and
  W.~Jones, {\it Handbook of Continued Fractions for Special Functions\/},
  Springer Netherlands (2008).
\bibitem{Jac}[Jac] L.~Jacobsen (ed.), {\it Analytic theory of continued
  fractions III\/}, Lecture Notes in Math.~{\bf 1406}, Springer (1989).
\bibitem{Jea}[Jea] R.~Andr\'e-Jeannin, {\it Irrationalit\'e de la somme
  des inverses de certaines suites r\'ecurrentes\/},
  C.~R.~Acad.~Sci.~{\bf 308} (1989), 539--541.
\bibitem{Jon-Thr}[Jon-Thr] W.~Jones and W.~Thron, {\it Continued fractions.
  Analytic theory and applications\/}, Enc.~Math.~and Appl., Addison--Wesley
  (1980).
\bibitem{JTW}[JTW] W.~Jones, W.~Thron, and H.~Waadeland (eds), {\it Analytic
  theory of continued fractions\/}, Lecture Notes in Math.~{\bf 932},
  Springer (1982).
\bibitem{Kho}[Kho] A.~Khovanskii, {\it The application of continued fractions\/}, Noordhoff (1963).
\bibitem{Khr}[Khr] S.~Khrushchev, {\it Orthogonal polynomials and continued
  fractions, From Euler's point of view\/}, Encyclopedia Math.~and Appl.~{\bf 122}, Cambridge (2008).
\bibitem{Lor-Waa}[Lor-Waa] L.~Lorentzen and H.~Waadeland, {\it Continued
  fractions\/}, Atlantis Studies in Math.~{\bf 1} (2008).
\bibitem{PWZ}[PWZ] M.~Petkovsek, H.~Wilf, and D.~Zeilberger, {\it $A=B$\/},
  A.~K.~Peters (1996).
\bibitem{Sti}[Sti] T.~Stieltjes, {\it Recherches sur les fractions continues\/}, Ann.~Fac.~Sci.~Toulouse {\bf 8} (1894), J1--122; {\bf 9} (1895), A1--47.
\bibitem{Thr}[Thr] W.~Thron (ed.), {\it  Analytic theory of continued
  fractions II\/}, Lecture Notes in Math.~{\bf 1199}, Springer (1986).
\bibitem{Wal}[Wal] H.~Wall, {\it Analytic theory of continued fractions\/},
  van Nostrand (1948).
\end{thebibliography}
\end{document}